\documentclass[11pt]{article}
\usepackage{amsmath,amsfonts,amscd,bezier}
\usepackage{latexsym,amssymb}
\usepackage[latin1]{inputenc}
\usepackage[english]{babel}

\newtheorem{theorem}{Theorem}[section]
\newtheorem{lemma}[theorem]{Lemma}
\newtheorem{proposition}[theorem]{Proposition}
\newtheorem{corollary}[theorem]{Corollary}
\newtheorem{definition}[theorem]{Definition}
\newtheorem{example}[theorem]{Example}
\newtheorem{remark}[theorem]{Remark}
\newcommand{\Dem}{\noindent{\it Proof.}\ \ }
\newcommand{\cqd}{{\hfill $\rule{2mm}{2mm}$}\vspace{.5cm}}

\begin{document}

\title{$\omega$-Symplectic algebra and Hamiltonian vector fields}

%
%
%

\author{Patr\'{\i}cia~H.~Baptistelli, Maria Elenice R. Hernandes\\
{\small Departamento de Matem\'atica, Universidade de Maring\'a}\\
{\small Av. Colombo, 5790, 87020-900, Maring\'a - PR, Brazil \footnote{Email address: phbaptistelli@uem.br (corresponding author), merhernandes@uem.br}}
\and Eralcilene Moreira Terezio\\
{\small Instituto Latino-Americano de Ci\^encias da Vida e da Natureza}\\ 
{\small Universidade Federal da Integra\c{c}\~ao Latino-Americana} \\ {\small Av. Tancredo Neves, 6731,  85867-900, Foz do Igua\c cu - PR, Brazil \footnote{Email address: eralcilene.terezio@unila.edu.br}}}

\date{}
\maketitle

\begin{abstract}
The purpose of this paper is presenting an algebraic theoretical basis for the study of $\omega$-Hamiltonian vector fields. We introduce the concepts of $\omega$-symplectic group and $\omega$-semisymplectic group, and describe some of their properties. We show that the Lie algebra of such groups is a useful tool in the recognition and construction of $\omega$-Hamiltonian vector fields defined on a symplectic vector space $(V,\omega)$ with respect to coordinates that are not necessarily symplectic. 

\end{abstract}



\section{Introduction}
An important tool in the study of symplectic vector spaces  is the well-known {symplectic group}  $Sp(n;\mathbb{R})$, that is, the set of all matrices $B\in{\mathbb{M}}_{2n}(\mathbb{R})$ such that $B^TJB= J$, where \begin{equation}\label{Jnporn}
	J=\left[\begin{array}{cc}0&I_n\\
		-I_n&0\end{array}\right],
\end{equation} with $I_n$ the identity matrix of order $n$. The group $Sp(n;\mathbb{R})$ has been studied by many authors and several of its properties can be found in \cite{NJ, Zehnder}. A symplectic vector space is a pair $(V,\omega)$, where $V$ is a finite-dimensional real vector space and $\omega:V\times V\to\mathbb{R}$ is an alternating and non-degenerate bilinear form. In this case, the dimension of $V$ is even and the matrix $[\omega]$ of $\omega$ relative to a basis of $V$ is a skew-symmetric matrix satisfying $\det[\omega]\neq0$. Any symplectic vector space $(V,\omega)$ admits a basis $\mathcal{C}$, called symplectic basis, such that $[\omega]_{\mathcal{C}} = J$. If $(V_1,\omega_1)$ and $(V_2,\omega_2)$ are symplectic vector spaces, then the set of linear maps $\xi:V_1\to V_2$ preserving the bilinear forms, that is, $$\omega_2(\xi(\mathbf{u}),\xi(\mathbf{v})) = \omega_1(\mathbf{u},\mathbf{v})$$ for all $\mathbf{u}, \mathbf{v}\in V_1$, is an interesting class of maps, called symplectic maps. When $\xi$ is an isomorphism, we say that $\xi$ is a symplectomorphism. A classical result in this context is that any symplectic vector space $(V,\omega)$ of dimension $2n$ is symplectomorphic to $(\mathbb{R}^{2n},\omega_0)$, where
\[\omega_0(\mathbf{u},\mathbf{v}) = \displaystyle\sum_{i=1}^n(x_iy_{n+i}-x_{n+i}y_i),\] with $\mathbf{u}=(x_1,\ldots,x_{2n})$ and $\mathbf{v}=(y_1,\ldots,y_{2n})$ in $\mathbb{R}^{2n}$ (see \cite{McDuff, Cannas, Zehnder}). In this case, the elements of the symplectic group $Sp(n;\mathbb{R})$ are matrices of symplectic operators in $(\mathbb{R}^{2n},\omega_0)$ and vice versa. 

This paper deals with the study of the symplectic algebra as a support to obtain results related to Hamiltonian vector fields defined on an arbitrary symplectic vector space $(V,\omega)$. An $\omega$-Hamiltonian vector field  is a map $X_H: V\to V$ given by\[X_H(x) = ([\omega]^{-1})^T \nabla H(x),\]for all $x \in V$, where $H:V\to\mathbb{R}$ is a smooth function. Our main motivation is the reference \cite{Terezio} in which we present normal forms for $\omega$-Hamiltonian vector fields on $(V,\omega)$ under the action of a group of symmetries and reversing symmetries. In this context, it was interesting to consider the matrix $[\omega]$ with respect to a not necessarily symplectic basis, since the symmetries and reversing symmetries of vector fields may not be preserved by symplectomorphisms. In \cite{buzzi1, Ricardo2013, Ricardo2011, Robinson, Sev} the authors also study Hamiltonian vector fields for other classes of matrices different than $J$.

In this work we consider an arbitrary symplectic vector space $(V,\omega)$. We show that if $[\omega]^2 = -I_{2n}$, then $(V,\omega)$ reflects similar properties to those of the classic symplectic space $(\mathbb{R}^{2n},\omega_0)$, but not in general. We introduce the concepts of $\omega$-symplectic group and $\omega$-semisymplectic group by describing their Lie algebra and some of their properties. The $\omega$-symplectic group $Sp_{\omega}(n;\mathbb{R})$ generalizes the group $Sp(n;\mathbb{R})$ and its matrices are in correspondence with linear symplectic maps. Similarly, we can consider a local diffeomorphism that reverses the symplectic form of the vector space, called antisymplectic map. We define the $\omega$-semisymplectic group
$\Omega_n = Sp_{\omega}(n;\mathbb{R}) \ \dot{\cup} \  Sp_{\omega}^{-1}(n;\mathbb{R})$, where $Sp^{-1}_{\omega}(n;\mathbb{R})$ is the set of antisymplectic matrices, which provides all algebraic information ne\-cessary for the study of symplectic and antisymplectic linear operators. These maps arise naturally in physical systems and are closely related to semisymplectic actions, which appear in several studies involving symmetric Hamiltonian systems, as for instance in \cite{Alomair, Terezio, buono, buzzi1, Montaldi2000, Montaldi, Claudia}.

%

We prove that the elements of the Lie algebra $\mathfrak{sp}_{\omega}(n;{\mathbb{R}})$ of the $\omega$-semisymplectic group are related to $\omega$-Hamiltonian vector fields. In particular, Proposition \ref{relacaomatrizcampo3} states that if $X$ is an $\omega$-Hamiltonian vector field then its linearization $dX_0$ belongs to ${\mathfrak{sp}}_{\omega}(n;\mathbb{R})$. Moreover, we characterize $X$ in terms of this Lie algebra when $X$ is linear and polynomial (Proposition \ref{relacaomatrizcampo2} and Theorem \ref{relacaopolinomiais}, respectively). The construction done in Proposition \ref{relacaomatrizcampo2} is used in Corollary \ref{relacaomatrizcampo2barra} to obtain families of nonlinear $\omega$-Hamiltonian vector fields whose linear part belongs to ${\mathfrak{sp}}_{\omega}(n;\mathbb{R})$.

This paper is organized as follows: in Section \ref{secaoaplicacoessimpleticas}, we present a brief study of the $\lambda$-symplectic maps. In Subsection \ref{subsec:Symplectic group}, we present the $\omega$-symplectic group $Sp_{\omega}(n;\mathbb{R})$, its Lie algebra ${\mathfrak{sp}}_{\omega}(n;\mathbb{R})$ and properties of the elements of both sets. In Subsection \ref{sec:semisymplectic} we introduce the $\omega$-semisymplectic group $\Omega_n$ and describe some of its algebraic and topological properties. In Section \ref{secaocamposmatrizes}, we relate $\omega$-Hamiltonian vector fields to elements of ${\mathfrak{sp}}_{\omega}(n;\mathbb{R})$.

\section{$\lambda$-Symplectic maps}\label{secaoaplicacoessimpleticas}
Let $(V,\omega)$ be a symplectic real vector space of dimension $2n$. In this section, we define $\lambda$-symplectic maps on $(V,\omega)$ and present some of their properties. The most interesting cases in the study of Hamiltonian vector fields are $\lambda = 1$ (symplectic maps) and $\lambda = -1$ (antisymplectic maps). The results of this section can naturally be extended to a symplectic manifold.

In what follows, $GL(2n)$ denotes the group of invertible matrices of order $2n$.  

\begin{definition}\label{conjlamdasimpletico}
	Given a skew-symmetric matrix $A \in GL(2n)$ and $\lambda\in\mathbb{R}^*$, we define the set\[{Sp}_{A}^{\lambda}(n;{\mathbb{R}})=\left\{B\in {\mathbb{M}}_{2n}(\mathbb{R}): B^TAB  = \lambda A\right\}.\]
\end{definition}

It is possible to prove that $Sp_{A}^{\lambda}(1;\mathbb{R})$ is the set of $2\times 2$ matrices with determinant equal to $\lambda$. In particular, $Sp_{A}^{1}(1;\mathbb{R}) = Sl(2;\mathbb{R})$, where $Sl(2n;\mathbb{R})$ is the well-known set of matrices of order $2n$ with determinant equal to 1. However, in general, $Sp_{A}^{1}(n;\mathbb{R})$ does not coincide with the set $Sl(2n;\mathbb{R})$. Indeed, the matrix
\[
B = \left[\begin{array}{cccc}
	1&0&0&0\\
	0&-2&0&0\\
	0&0&\frac{1}{2}&0\\
	0 &0&0&-1
\end{array}
\right]\] belongs to $Sl(4,\mathbb{R})$, but $B\notin Sp_{A}^1(2;\mathbb{R})$, for any skew-symmetric and invertible matrix $A$.

Linear operators on ${\mathbb{R}}^{2n}$ which preserve $\omega_0$ are in correspondence with elements of the symplectic group. A natural question that arises is: under what conditions does the matrix of a linear operator belong to ${Sp}_{A}^{\lambda}(n;{\mathbb{R}})$? In order to answer this question, we present the following definition:
\begin{definition}\label{definicaolambdasimpletica} Let $(V,\omega)$ be a symplectic vector space. Given $\lambda\in\mathbb{R}^*$, a differentiable map $\xi:V\to V$ is called $\lambda$-symplectic if\begin{center}
		$\xi^*\omega = \lambda\omega$,
	\end{center}where $\xi^*\omega$ is the pullback of $\omega$ by $\xi$. The map $\xi$ is called symplectic if $\lambda =1$ and antisymplectic if $\lambda =-1$.
\end{definition}

Note that any $\lambda$-symplectic linear operator is invertible. Consequently, a $\lambda$-symplectic map is an immersion and, therefore, its local inverse exists. 

\begin{proposition}\label{matrizaplicacaosimpletica} Let ${\mathcal{C}}$ be a basis of $(V,\omega)$. If $\xi:V\to V$ is a $\lambda$-symplectic linear operator, then $[\xi]_{\mathcal{C}}\in {Sp}_{[\omega]_{\mathcal{C}}}^{\lambda}(n;{\mathbb{R}})$, where $[\xi]_{\mathcal{C}}$ and $[\omega]_{\mathcal{C}}$ denote the matrices of $\xi$ and of $\omega$ relative to the basis ${\mathcal{C}}$, respectively. Conversely, if $B\in{Sp}_{[\omega]_{\mathcal{C}}}^{\lambda}(n;{\mathbb{R}})$, then the linear operator $\xi:V\to V$ such that $[\xi]_{\mathcal{C}} = B$ is $\lambda$-symplectic.
\end{proposition}

\Dem If $\xi$ is a $\lambda$-symplectic linear operator, then for all $\mathbf{u},\mathbf{v}\in V$ we have $\omega (\xi(\mathbf{u}),\xi(\mathbf{v})) = \lambda\omega (\mathbf{u},\mathbf{v})$. In matrix terms, \begin{center}
		$\mathbf{u}^T[\xi]_{\mathcal{C}}^T[\omega]_{\mathcal{C}}[\xi]_{\mathcal{C}}\mathbf{v}=([\xi]_{\mathcal{C}}\mathbf{u})^T[\omega]_{\mathcal{C}}[\xi]_{\mathcal{C}}\mathbf{v} = \omega (\xi(\mathbf{u}),\xi(\mathbf{v})) = \lambda\omega (\mathbf{u},\mathbf{v})= \mathbf{u}^T\lambda[\omega]_{\mathcal{C}} \mathbf{v}$.
	\end{center}Hence $[\xi]_{\mathcal{C}}^T[\omega]_{\mathcal{C}} [\xi]_{\mathcal{C}} = \lambda[\omega]_{\mathcal{C}} $. On the other hand, given $B\in \mathbb{M}_{2n}(\mathbb{R})$ such that $B^T[\omega]_{\mathcal{C}}B = \lambda[\omega]_{\mathcal{C}}$, by considering the linear operator $\xi:V\to V$ such that $[\xi]_{\mathcal{C}} = B$, we have\[\omega (\xi(\mathbf{u}),\xi(\mathbf{v})) = (B\mathbf{u})^T[\omega]_{\mathcal{C}}B\mathbf{v} = \mathbf{u}^TB^T[\omega]_{\mathcal{C}}B\mathbf{v} = \mathbf{u}^T\lambda[\omega]_{\mathcal{C}}\mathbf{v} =  \lambda\omega (\mathbf{u},\mathbf{v}).\] Thus, $\xi$ is $\lambda$-symplectic. \cqd

If $[\omega]$ denotes the matrix of $\omega$ relative to a basis $\mathcal{B}$, in order to simplify the notation we denote the set ${Sp}_{[\omega]}^{\lambda}(n;{\mathbb{R}})$ simply by\begin{equation}\label{defigruposimpleticocorrecao}
	{Sp}_{\omega}^{\lambda}(n;{\mathbb{R}}) =\left\{B\in {\mathbb{M}}_{2n}(\mathbb{R}): B^T[\omega]B  = \lambda [\omega]\right\}.
\end{equation} The matrices in ${Sp}_{\omega}^{-1}(n;{\mathbb{R}})$ are called \emph{$\omega$-antisymplectic}. In addition, we denote ${Sp}_{\omega}^{1}(n;{\mathbb{R}})$ by ${Sp}_{\omega}(n;{\mathbb{R}})$ and their matrices are called \emph{$\omega$-symplectic}. 

One property of the matrices in ${Sp}_{\omega}^{\lambda}(n;{\mathbb{R}})$ is the following:
\begin{proposition}\label{detanti}
	Let $(V,\omega)$ be a symplectic vector space. If $B\in{Sp}_{\omega}^{\lambda}(n;{\mathbb{R}})$, then $\det B = \lambda^n$.
\end{proposition}

\Dem By Proposition \ref{matrizaplicacaosimpletica}, if $B\in{Sp}_{\omega}^{\lambda}(n;{\mathbb{R}})$, then there exists a $\lambda$-symplectic linear operator $\xi:V\to V$ such that $[\xi]_{\mathcal{B}} = B$. Consider the constant $2n$-form $\varphi = \omega^n$ defined as the exterior product of $n$ factors of $\omega$. We obtain
	\begin{center}
		$\xi^*\varphi = \xi^*\omega\wedge\ldots\wedge \xi^*\omega=\lambda\omega\wedge\ldots\wedge\lambda\omega=\lambda^n\varphi$.
	\end{center}  
	Since $\xi^*\varphi  = (\det B)\varphi$ (see \cite[Proposition 14.20]{lee}), we have $\det B = \lambda^n$. \cqd 

In particular, given an $\omega$-antisymplectic matrix $B\in \mathbb{M}_{2n}(\mathbb{R})$, then $\det B = 1$ if $n$ is even and $\det B = -1$ if $n$ is odd. In addition, given an $\omega$-symplectic matrix $B\in \mathbb{M}_{2n}(\mathbb{R})$, $\det B = 1$ for all $n$. In this sense, the Proposition \ref{detanti} is a generalization of the known Liouville's theorem. 

In the next results we explore some properties of $\lambda$-symplectic maps.
\begin{proposition}\label{deltaSpnA}Let $\delta:V\to V$ be a $\lambda_1$-symplectic map. A map $\xi:V\to V$ is $\lambda_2$-symplectic if and only if $\delta\circ \xi$ is $\lambda_1\lambda_2$-symplectic.
\end{proposition}

\Dem By hypothesis $\delta^*\omega = \lambda_1\omega$. If $\xi$ is a $\lambda_2$-symplectic map, that is, $\xi^*\omega =\lambda_2\omega $, then
	\[(\delta\circ \xi)^*\omega  = \xi^*({\delta}^*\omega ) =\xi^*(\lambda_1\omega ) = \lambda_1(\xi^*\omega )= \lambda_1(\lambda_2\omega ) = \lambda_1\lambda_2\omega.\]Conversely, if $\delta\circ \xi$ is a $\lambda_1\lambda_2$-symplectic map, we have\[\lambda_1\lambda_2\omega  = (\delta\circ \xi)^*\omega  = \xi^*(\delta^*\omega ) = \xi^*(\lambda_1\omega ) = \lambda_1(\xi^*\omega ),\]
	that is, $\xi^*\omega  = \lambda_2\omega$.\cqd 

\begin{corollary}\label{CorolariosEnumerados}
	\begin{itemize}
			\item[(i)] A map $\delta:V\to V$ is antisymplectic if and only if $\delta\circ \xi$ is antisymplectic for any symplectic map $\xi:V\to V$.
		\item[(ii)] A map $\delta:V\to V$ is $\lambda$-symplectic if and only if its (local) inverse $\delta^{-1}$ is a $\lambda^{-1}$-symplectic map.
		\item[(iii)] Given an antisymplectic map $\delta:V\to V$, then $\delta^k$ is  symplectic if $k$ is even and antisymplectic if $k$ is odd.
		\end{itemize}
\end{corollary}

\section{$\omega$-Symplectic algebra}\label{secaoalgebrasimpletica}

In this section, we present a generalization of the symplectic group and some of its algebraic properties. For references of these results on $(\mathbb{R}^{2n},\omega_0)$ see for instance \cite{Hofer, McDuff, Robinson, Zehnder}. We also define the $\omega$-semisymplectic group and determine its Lie algebra, which is an interesting tool in the study of $\omega$-Hamiltonian vector fields (see Section \ref{secaocamposmatrizes}).

\subsection{$\omega$-Symplectic group and its Lie algebra} \label{subsec:Symplectic group}

We start showing that the set $Sp_{\omega}(n;\mathbb{R})$ defined in (\ref{defigruposimpleticocorrecao}) for $\lambda = 1$ is a linear Lie group, that is, a closed subgroup of $GL(2n)$. Note that ${Sp}_{\omega}(n;{\mathbb{R}})$ is a subgroup of $GL(2n)$ whose identity element is $I_{2n}$. Moreover $Sp_{\omega}(n;\mathbb{R})$ is closed, since $Sp_{\omega}(n;\mathbb{R}) = f^{-1}(\{[\omega]\})$ where $f:GL(2n)\to GL(2n)$ is the continuous map given by $f(A) = A^T[\omega]A$. The Lie group \begin{center}
	$Sp_{\omega}(n;\mathbb{R})=\left\{B\in GL(2n): B^T[\omega] B = [\omega] \right\}$\index{$Sp_{\omega}(n;\mathbb{R})$}
\end{center}is called \emph{$\omega$-symplectic group}.

\begin{example}\label{exacaosimpleticaso2}{\rm
		Let $({\mathbb{R}}^2,\omega)$ be an arbitrary symplectic vector space. We consider the special orthogonal group $SO(2)$ generated by the rotation matrices  \begin{equation}\label{NUMERE}
			R_{\theta} = \left[\begin{array}{cr}
				\cos\theta&-\sin\theta\\
				\sin\theta&\cos\theta
			\end{array}\right].
		\end{equation} Since $\det R_{\theta}=1$ for all $\theta \in [0,2\pi)$, it follows that $SO(2)$ is a subgroup of $Sp_{\omega}(1,{\mathbb{R}})$.

}\end{example}

Note that if we choose a symplectic basis of $V$, then $B\in Sp_{\omega}(n;{\mathbb{R}})$ if and only if $B^T\in Sp_{\omega}(n;{\mathbb{R}})$. More generally, we have:

\begin{proposition} \label{prop:description} Let $(V,\omega)$ be a symplectic vector space. If $[\omega]^2= -I_{2n}$, then\begin{center}
		$Sp_{\omega}(n;\mathbb{R})=\left\{B\in GL(2n): B[\omega] B^T=B^T[\omega] B = [\omega] \right\}$.
	\end{center}
\end{proposition}

\Dem If $[\omega]^2= -I_{2n}$, then $[\omega]^{-1} = -[\omega]$. We will show that if $B\in Sp_{\omega}(n;\mathbb{R})$, then $B^T\in Sp_{\omega}(n;\mathbb{R})$. Indeed, if $B\in Sp_{\omega}(n;\mathbb{R})$, then $B^{-1}\in Sp_{\omega}(n;\mathbb{R})$. Hence
	$$	-[\omega]  = [\omega]^{-1} = ((B^{-1})^T[\omega]B^{-1})^{-1} = B(-[\omega])((B^{-1})^T)^{-1} = -B[\omega]B^T, 
	$$ that is, $[\omega]  = B[\omega]B^T$, for all $B\in Sp_{\omega}(n;{\mathbb{R}})$.\cqd 

Since the symplectic group $Sp(n;\mathbb{R})$ is the $\omega$-symplectic group in a symplectic basis of
$(V,\omega)$, by Proposition \ref{prop:description} we have $$
Sp(n;\mathbb{R}) = \left\{B\in GL(2n): BJ B^T=B^TJ B = J \right\},$$ where $J$ is given in (\ref{Jnporn}).  The Lie groups $Sp_{\omega}(n;\mathbb{R})$ and $Sp(n;\mathbb{R})$ are isomorphic, since the map $f:Sp_{\omega}(n;\mathbb{R})\to Sp(n;\mathbb{R})$ given by $f(B) = P^{-1}BP$ is an isomorphism of Lie groups, where $P$ is the change of basis matrix.

\begin{remark}\label{obsgruposisomorfos}\textrm{For $\lambda\in{\mathbb{R}}^*$ fixed, the bijective map $f:Sp_{\omega}^{\lambda}(n;\mathbb{R})\to Sp_J^{\lambda}(n;\mathbb{R})$ given by $f(B) = P^{-1}BP$ as above  guarantees that $Sp_{\omega}^{-1}(n;\mathbb{R})$ is a non-empty set since, for instance, the matrix \[\left[\begin{array}{cc}
			0&I_n\\
			I_n&0
		\end{array}\right]\]belongs to $Sp_{J}^{-1}(n;\mathbb{R})$.
}\end{remark}

\begin{proposition}\label{provaalgLie}The Lie algebra ${\mathfrak{sp}}_{\omega}(n;{\mathbb{R}})$ of the $\omega$-symplectic group satisfies\begin{align}\label{eqAlgebraspomega}
		{\mathfrak{sp}}_{\omega}(n;{\mathbb{R}}) & =\left\{L\in{\mathbb{M}}_{2n}({\mathbb{R}}): L^T[\omega]  +[\omega] L=0\right\}\nonumber\\
		& = \left\{L\in{\mathbb{M}}_{2n}({\mathbb{R}}): ([\omega]^{-1})^{T}L^T[\omega]  = L\right\}.
	\end{align}
\end{proposition}

\Dem The second equality in (\ref{eqAlgebraspomega}) follows from the property that $[\omega]$ is a skew-symmetric matrix. By definition of Lie algebra we have that\begin{equation*}
		{\mathfrak{sp}}_{\omega}(n;{\mathbb{R}}) = \left\{L\in {\mathbb{M}}_{2n}({\mathbb{R}}): \exp({tL})\in Sp_{\omega}(n;\mathbb{R}), \ \forall t\in{\mathbb{R}}\right\}.
	\end{equation*}Hence, if  $L\in{\mathfrak{sp}}_{\omega}(n;{\mathbb{R}})$, then $(\exp(tL))^T[\omega]\exp(tL) = [\omega]$  for all $t\in{\mathbb{R}}$, that is, $[\omega]  = \exp(tL^T)[\omega] \exp(tL)$. Thus
	\[0  = \dfrac{\textrm{d}}{{\textrm{d}}t}\,\left(\exp(tL^T)[\omega] \exp(tL)\right)\bigg|_{t=0}
	= L^T[\omega] +[\omega] L.\]Conversely, if $L\in{\mathbb{M}}_{2n}({\mathbb{R}})$ satisfies $([\omega]^{-1})^{T}L^T[\omega] = L$, then $([\omega] ^{-1})^T(tL^T)[\omega] = tL$ for each $t\in{\mathbb{R}}$. In this way, \[\begin{array}{rcl}
		\exp(tL^T)[\omega]\exp(tL)&=&\exp(tL^T)[\omega]\exp(([\omega]^{-1})^T(tL^T)[\omega])\\
		&=& \exp(tL^T)[\omega]\exp({[\omega]^{-1}(-tL)^T[\omega]})\\
		&=& \exp(tL^T)[\omega][\omega]^{-1}\exp({(-tL)^T})[\omega] = [\omega],
	\end{array}\]that is, $L\in{\mathfrak{sp}}_{\omega}(n;{\mathbb{R}})$.\cqd  

In this work, ${\mathfrak{sp}}_{\omega}(n;{\mathbb{R}})$ is called $\omega$-symplectic algebra and its elements are called \emph{$\omega$-Hamiltonian matrices}. In the particular case in which $[\omega] = J$, we assume the classical notation ${\mathfrak{sp}}(n;{\mathbb{R}})$ for this Lie algebra.


\begin{remark}
	\label{remark:tranpose}
	Note that if $[\omega]^2 = -I_{2n}$, then\[{\mathfrak{sp}}_{\omega}(n;{\mathbb{R}}) = \{L\in{\mathbb{M}}_{2n}({\mathbb{R}}): [\omega] L^T[\omega]  = L\}.\]In this case, if $L\in{\mathfrak{sp}}_{\omega}(n;{\mathbb{R}})$, then $L^T\in{\mathfrak{sp}}_{\omega}(n;{\mathbb{R}})$. This property is not true in general. Indeed, given the matrices \begin{equation}\label{Snaosimpletico}
		L = \left[\begin{array}{cccc}
			-1&1&-1&2\\
			3&0&4&1\\
			-1&2&0&2\\
			3&1&1&1
		\end{array}\right] \quad  \textrm{ and} \quad [\omega]=\left[\begin{array}{cccc}
			0&1&0&2\\
			-1&0&-1&0\\
			0&1&0&1\\
			-2&0&-1&0
		\end{array}\right],
	\end{equation}it is easy to verify that $L\in{\mathfrak{sp}_{\omega}(2,\mathbb{R})}$, but $L^T\notin{\mathfrak{sp}_{\omega}(2,\mathbb{R})}$. \end{remark}

By a similar argument as we used for symplectic groups, we can prove that the Lie algebras ${\mathfrak{sp}}_{\omega}(n;{\mathbb{R}})$ and ${\mathfrak{sp}}(n;{\mathbb{R}})$ are isomorphic. It follows from this proof that $\omega$-Hamiltonian matrices and $\omega_0$-Hamiltonian matrices are conjugated. Moreover, the trace of an $\omega_0$-Hamiltonian matrix is zero (see \cite[Proposition 1.6]{Robinson}). Since the trace of a matrix is invariant by conjugation, we conclude the following:

\begin{proposition}
	\label{prop:trace}
	The trace of any $\omega$-Hamiltonian matrix is zero.
\end{proposition}

\subsection{$\omega$-Semisymplectic group and its Lie algebra}

\label{sec:semisymplectic}

In this subsection, we introduce the $\omega$-semisymplectic group and characterize its Lie algebra. Firstly, we observe that the set ${Sp}_{\omega}^{\lambda}(n;{\mathbb{R}})$ does not have a group structure for all $\lambda\neq 1$. However, Proposition \ref{deltaSpnA} and Corollary \ref{CorolariosEnumerados} $(ii)$ ensure that \[\displaystyle\Lambda = \bigcup_{\lambda\in\mathbb{R}^*}{Sp}_{\omega}^{\lambda}(n;{\mathbb{R}})\]is a subgroup of $GL(2n)$ and therefore its closure is a linear Lie group. 

By Remark \ref{obsgruposisomorfos}, the set of $\omega$-antisymplectic matrices is non-empty, so we consider $\delta \in {Sp}_{\omega}^{-1}(n;{\mathbb{R}})$ arbitrary and $\delta Sp_{\omega}(n;\mathbb{R}) = \left\{\delta B: B\in Sp_{\omega}(n;\mathbb{R})\right\}$. It follows from Proposition \ref{deltaSpnA} and Corollary \ref{CorolariosEnumerados} that\[\delta Sp_{\omega}(n;\mathbb{R}) = {Sp}_{\omega}^{-1}(n;{\mathbb{R}}) = \left\{B\in\mathbb{M}_{2n}(\mathbb{R}): B^T[\omega] B = -[\omega] \right\}.\] Since $[\omega]\neq0$, we define the subset of $\Lambda$ given by the disjoint union
\begin{equation}
	\label{eq:omegan}
	\Omega_n = Sp_{\omega}(n;\mathbb{R}) \ \dot{\cup} \  Sp_{\omega}^{-1}(n;\mathbb{R}).\end{equation}
As $Sp_{\omega}(n;\mathbb{R})$ is a group, from Corollary \ref{CorolariosEnumerados} we have that $\Omega_n$ is also a group. Moreover $\Omega_n = f^{-1}(\{[\omega] ,-[\omega] \})$, where $f:GL(2n)\to GL(2n)$ is the continuous map given by $f(B) = B^T[\omega] B$. Thus $\Omega_n$ is a linear Lie group, called \emph{$\omega$-semisymplectic group}. 

We can also see that $Sp_{\omega}(n;\mathbb{R})$ is a Lie subgroup of $\Omega_n$ of index $2$, because $Sp_{\omega}^{-1}(n;\mathbb{R}) = \delta Sp_{\omega}(n;\mathbb{R})$. Then there exists a group epimorphism $\sigma:\Omega_n\to\mathbb{Z}_2$ defined by
\[
\sigma(B) = \left\{
\begin{array}{rl}
	1, &\text{if} \ B\in Sp_{\omega}(n;\mathbb{R})\\
	-1, &\text{if} \ B\in Sp_{\omega}^{-1}(n;\mathbb{R})
\end{array}\right. .\]
According to Proposition \ref{detanti}, if $B\in Sp_{\omega}(n;\mathbb{R})$, then $\det B = 1$ and if $B\in Sp_{\omega}^{-1}(n;\mathbb{R})$, then $\det B = (-1)^n$. Thus if $n$ is odd, then $\sigma = \det$.


\begin{example}\label{exemplos Cap4}{\rm Consider the symplectic vector space $(\mathbb{R}^{2n}, \omega)$ and suppose that $B$ and $[\omega]$ are block matrices in $\mathbb{M}_{2n}(\mathbb{R})$ given by\[
		B= \left[\begin{array}{cccc}
			B_1&&&\\
			&B_2&&\\
			&&\ddots&\\
			&&&B_n
		\end{array}\right] \quad \text{and} \quad [\omega]= \left[\begin{array}{cccc}
			A_1&&&\\
			&A_2&&\\
			&&\ddots&\\
			&&&A_n
		\end{array}\right],\] with $A_j,B_j\in\mathbb{M}_{2}(\mathbb{R})$, where \begin{equation} \label{matrixAj} A_j = \left[\begin{array}{cc}
				0&a_j\\
				-a_j&0
			\end{array} \right],\end{equation} for all $j=1,\ldots,n$. We have that $B\in\Omega_n$ if and only if $B_j \in\Omega_1$,  for all $j=1,\ldots,n$. More specifically, $B\in Sp_{\omega}(n;\mathbb{R})$ if and only if $B_j \in Sp_{A_j}(1;\mathbb{R})$ and $B\in Sp_{\omega}^{-1}(n;\mathbb{R})$ if and only if $B_j \in Sp_{A_j}^{-1}(1;\mathbb{R})$, for all $j=1,\ldots,n$. In this case, $B\in Sp_{\omega}(n;\mathbb{R})$ if and only if $\det B_j = 1$, for all $j=1,\ldots,n$, and $B\in Sp_{\omega}^{-1}(n;\mathbb{R})$ if and only if $\det B_j = -1$, for all $j=1,\ldots,n$.}\end{example}

\begin{example}\label{exacaosimpleticao4}{\rm Let $({\mathbb{R}}^4,\omega)$ be a symplectic vector space, where the matrix of $\omega$ is given as in Example \ref{exemplos Cap4} for $n=2$. Consider the matrices \[\kappa =\left[\begin{array}{rc}
			-1&0\\
			0&1
		\end{array}\right], \quad \tau = \left[\begin{array}{cc}
			\kappa&0\\
			0&\kappa
		\end{array}\right] \quad \textrm{and}\quad
		\rho_{\theta} = \left[\begin{array}{cc}
			R_{\theta}&0\\
			0&R_{\theta}
		\end{array}\right],\]
		with $R_{\theta}$ given in (\ref{NUMERE}). Consider also the semidirect product $SO(2)\rtimes\mathbb{Z}_2^{\tau},$ where $\mathbb{Z}_2^{\tau}$ is the group generated by $\tau$ and $SO(2)$ is the subgroup of $SO(4)$ generated by matrices $\rho_{\theta}$. Since
		$R_{\theta}\in Sp_{A_j}(1;{\mathbb{R}})$ and $\kappa\in Sp_{A_j}^{-1}(1;{\mathbb{R}})$, for $A_j$ as in (\ref{matrixAj}) and $j=1,2$, then $\rho_{\theta}\in Sp_{\omega}(2;{\mathbb{R}})$ and $\tau\in Sp_{\omega}^{-1}(2;{\mathbb{R}})$. Therefore, $SO(2)\rtimes\mathbb{Z}_2^{\tau}$ is a subgroup of $\Omega_2.$
	}
\end{example}

\begin{lemma}\label{Spnaberto}The group $Sp_{\omega}(n;\mathbb{R})$ and the set $Sp_{\omega}^{-1}(n;\mathbb{R})$ are open and closed in $\Omega_n$.
\end{lemma}

\Dem Consider the continuous map $f:GL(2n)\to GL(2n)$ given by $f(B) = B^T[\omega]B$. Note that $Sp_{\omega}^{-1}(n;\mathbb{R}) = {f}^{-1}(\{-[\omega] \})$ and, therefore, $Sp_{\omega}^{-1}(n;\mathbb{R})$ is closed in $GL(2n)$. Thus, $GL(2n) - Sp_{\omega}^{-1}(n;\mathbb{R})$ is open in $GL(2n)$. Since $Sp_{\omega}(n;\mathbb{R}) = \Omega_n\cap\left(GL(2n) - Sp_{\omega}^{-1}(n;\mathbb{R})\right)$, we have that $Sp_{\omega}(n;\mathbb{R})$ is open in $\Omega_n$. Then $Sp_{\omega}^{-1}(n;\mathbb{R})$ is closed in $\Omega_n$. Similarly, we can conclude that $Sp_{\omega}^{-1}(n;\mathbb{R})$ is open in $\Omega_n$, which implies that the group $Sp_{\omega}(n;\mathbb{R})$ is also closed in $\Omega_n$.\cqd

An immediate consequence from Lemma \ref{Spnaberto} is that $\Omega_n$ is not connected. Moreover, $Sp_{\omega}(n;\mathbb{R})$ is the connected component of the identity of $\Omega_n$. This allows us to describe the Lie algebra  of $\Omega_n$, since the connected component of the identity of a Lie group $G$ have the same Lie algebra of $G$.

\begin{proposition}\label{mesmaalgebra}The Lie algebra of $\Omega_n$ is the set ${\mathfrak{sp}}_{\omega}(n;\mathbb{R})$ of the $\omega$-Hamiltonian matrices.
\end{proposition}


\section{Relation between $\omega$-Hamiltonian vector fields and ${\mathfrak{sp}}_{\omega}(n;\mathbb{R})$}\label{secaocamposmatrizes}


In this section, we relate $\omega$-Hamiltonian vector fields to elements of the Lie algebra ${\mathfrak{sp}}_{\omega}(n;\mathbb{R})$. Particularly, we show that ${\mathfrak{sp}}_{\omega}(n;\mathbb{R})$ simplifies the recognition of linear and polynomial $\omega$-Hamiltonian vector fields and it is a tool for the construction of nonlinear $\omega$-Hamiltonian vector fields. Some of
the results presented here are generalizations of \cite[Section 2.7]{Chow}, where the authors assume $[\omega] = J$ as in $(\ref{Jnporn})$. 

\begin{definition}\label{deficampoomegahamil}
	The $\omega$-Hamiltonian vector field associated with a smooth function $H:V\to\mathbb{R}$ is the unique vector field $X_H:V\to V$ such that $\omega(X_H(x),\cdot) = dH_x$, for all $x\in V$. In this case, $X_H$ can be written as
	\begin{equation}
		\label{deficampoHmatriz*}
		X_H(x) = ([\omega]^{-1})^{T}\nabla H(x)
	\end{equation}for all $x\in V$. The function $H$ is called a Hamiltonian function.
\end{definition}

We note that the uniqueness of $X_H$ is characterized by the non-degeneracy of $\omega$. It is possible to prove that the flows of a Hamiltonian vector field preserve the symplectic form (see \cite[Proposition V.24]{Zehnder}). In addition, according to \cite[Proposition V.26]{Zehnder}, if $\xi:V\to V$ is a symplectic map and $X_H$ is a Hamiltonian vector field, then the pushforward of $X_H$ by $\xi$ is also a Hamiltonian vector field. Now we present a natural generalization of this result.

\begin{theorem}\label{mudancacompostaartigo}
	Let $(V,\omega)$ be a symplectic vector space and $X_H:V\to V$ an $\omega$-Hamiltonian vector field associated with $H:V\to\mathbb{R}$. If $\xi:V\to V$ is a $\lambda$-symplectic map, then the induced vector field $\xi_*X_H$ of $X_H$ by $\xi$ satisfies\[\xi_*X_H = \lambda X_{H\circ \xi^{-1}},\]where $\xi^{-1}$ denotes the local inverse of $\xi$.
\end{theorem}

\Dem By Proposition V.11 and Proposition V.25 of \cite{Zehnder} it follows, respectively, that $d(\xi^*\omega) = \xi^*(d\omega)$ and $\xi^*(\omega(X_H(x), \cdot)) = (\xi^*\omega)((\xi^{-1})_*X_H(x), \cdot)$ for all $x \in V$. By Definition \ref{deficampoomegahamil}, we obtain \begin{align*}
		\omega({X_{H\circ \xi^{-1}}}(x),\cdot) = d(H\circ \xi^{-1})_x &= (\xi^{-1})^*(dH_x)\\
		& = (\xi^{-1})^*(\omega(X_H(x), \cdot))\\
		& = ((\xi^{-1})^*\omega)(\xi_*X_H(x),\cdot).
	\end{align*}
	By Corollary \ref{CorolariosEnumerados} $(ii)$, $\xi^{-1}$ is a $\lambda^{-1}$-symplectic map. So \begin{center}
		$ 0 = \omega({X_{H\circ \xi^{-1}}}(x),\cdot) - (\lambda^{-1}\omega)(\xi_*X_H(x),\cdot) = \omega(X_{H\circ \xi^{-1}}(x)-\lambda^{-1}\xi_*X_H(x),\cdot)$ \end{center} for all $x \in V$. Since $\omega$ is non-degenerate, we conclude the result.\cqd  

For what follows, we denote by $dX_x$ both the linear operator and its matrix relative to a basis of $(V,\omega)$.

\begin{proposition}\label{relacaomatrizcampo1}Let $X:V\to V$ be an $\omega$-Hamiltonian vector field. Then there exists a symmetric matrix $B$ satisfying $dX_0 = ([\omega]^{-1})^TB$. Conversely, given a symmetric matrix  $B$, there exists an $\omega$-Hamiltonian linear vector field $X:V\to V$ such that $dX_0 = ([\omega]^{-1})^TB$.
\end{proposition}

\Dem Since $X$ is an $\omega$-Hamiltonian vector field, by (\ref{deficampoHmatriz*}) there exists a function $H:V\to{\mathbb{R}}$ such that $X(x) = ([\omega]^{-1})^T\nabla H(x)$. Now, it is enough to consider $B = d^2H_0$, where $d^2H_0$ denotes the Hessian matrix of $H$ in $x=0$.
	
	Reciprocally, given $B = [b_{ij}]\in\mathbb{M}_{2n}(\mathbb{R})$ a symmetric matrix, we define the quadratic form $H:V\to{\mathbb{R}}$ as \[H(x) = \dfrac{1}{2}\mathbf{x}^TB\mathbf{x} = \dfrac{1}{2}\displaystyle\sum_{i,j=1}^{2n}x_ib_{ij}x_j,\]where $\mathbf{x}^T = [x_1 \ \ldots \ x_{2n}]$. The $\omega$-Hamiltonian vector field $X$ associated with $H$ satisfies 
	$$X(x) = ([\omega]^{-1})^{T}\nabla H(x) = ([\omega]^{-1})^{T}B\mathbf{x}$$ and therefore $dX_0 = ([\omega]^{-1})^TB$.\cqd  

As a direct consequence of the next result, it is possible to recognize if a vector field is not $\omega$-Hamiltonian only
by analyzing the trace of its linearization.

\begin{proposition}\label{relacaomatrizcampo3}If  $X$ is an $\omega$-Hamiltonian vector field, then $dX_x \in {\mathfrak{sp}}_{\omega}(n;{\mathbb{R}})$, for all $x\in V$.
\end{proposition}

\Dem If $X$ is an $\omega$-Hamiltonian vector field, then there exists a function $H:V\to {\mathbb{R}}$ such that $X(x) = ([\omega]^{-1})^{T}\nabla H(x)$, whence $dX_x = ([\omega]^{-1})^Td^2H_x$. Thus,\begin{center}
		$\begin{array}{rcl} (dX_x)^T[\omega]+[\omega] dX_x &=& (([\omega]^{-1})^{T}d^2H_x)^T[\omega]+[\omega]([\omega]^{-1})^{T}d^2H_x\\
			& =& d^2H_x[\omega]^{-1}[\omega]+[\omega](-[\omega])^{-1}d^2H_x = 0.
		\end{array}$
	\end{center}Hence, $dX_x\in\mathfrak{sp}_{\omega}(n;{\mathbb{R}})$ for all $x\in V$.  \cqd  

\begin{example}\label{naohamiltoniano2}{\rm
		By Propositions \ref{prop:trace} and \ref{relacaomatrizcampo3}, the vector field $X:\mathbb{R}^6\to\mathbb{R}^6$ given by \begin{align*}
			X(x) = &\left(x_1+x_3+x_5+x_1^2, x_2+x_4+x_6+x_2x_4x_6,\right.\\
			& \left.-x_3+x_5+x_3x_5, 3x_4, -3x_5+x_6-x_6^3, 2x_6\right),
		\end{align*}for all $x = (x_1,\ldots,x_6) \in \mathbb{R}^6$, is not an $\omega$-Hamiltonian vector field whatever the symplectic structure $\omega$ on $\mathbb{R}^6$, since the trace of the linearization $L=dX_0$ of $X$ is non-zero.}\end{example}

In \cite[Chapter 2, Lemma 6.2]{Chow}, the authors show that a linear vector field $X(x) = L\mathbf{x}$ on $({\mathbb{R}}^{2n},\omega_0)$ is Hamiltonian if and only if $L\in{\mathfrak{sp}}(n;{\mathbb{R}})$. The next result extends this equivalence for linear vector fields on an arbitrary symplectic vector space $(V,\omega)$.

\begin{proposition}\label{relacaomatrizcampo2}A matrix $L$ belongs to ${\mathfrak{sp}}_{\omega}(n;{\mathbb{R}})$ if and only if the linear vector field $X(x) = L\mathbf{x}$ is $\omega$-Hamiltonian.
\end{proposition}

\Dem Given $L\in{\mathfrak{sp}}_{\omega}(n;{\mathbb{R}})$, we define the symmetric matrix $B := [{\omega}]^{T}L$ and the function $H(x) = \frac{1}{2}\mathbf{x}^TB\mathbf{x}$. We have that the $\omega$-Hamiltonian vector field associated with $H$ satisfies \[
	X(x) = ({[\omega]}^{-1})^{T}B\mathbf{x} = ({[\omega]}^{-1})^{T}{[\omega]}^{T}L\mathbf{x} = L\mathbf{x}.\]The sufficient condition is immediate from Proposition \ref{relacaomatrizcampo3}.\cqd  

It is easy to see that if $B\in\mathbb{M}_{2n}(\mathbb{R})$ is a symmetric matrix, then $([\omega]^{-1})^TB \in {\mathfrak{sp}}_{\omega}(n;{\mathbb{R}})$. By proof of the previous proposition we obtain another characterization for the Lie algebra ${\mathfrak{sp}}_{\omega}(n;{\mathbb{R}})$, namely\begin{center}
	${\mathfrak{sp}}_{\omega}(n;{\mathbb{R}}) = \left\{([\omega]^{-1})^TB: \ B\in{\mathbb{M}}_{2n}({\mathbb{R}}) \ \ \textrm{is} \ \ \textrm{symmetric}\right\}$.
\end{center}

We show in the next result that the construction done in Proposition \ref{relacaomatrizcampo2} can be used to obtain families of nonlinear $\omega$-Hamiltonian vector fields. For this, we introduce the following definition. Given a function $F:\mathbb{R}^m\to\mathbb{R}$ of class $C^{\infty}$ and $a\in\mathbb{R}^m$, we denote by $j^kF(a)$ the Taylor polynomial of $F(x+a)-F(a)$ at $0$ of order $k$, called \emph{$k$-jet of $F$ at $a$}.

\begin{corollary}\label{relacaomatrizcampo2barra} If $L\in{\mathfrak{sp}}_{\omega}(n;{\mathbb{R}})$, then there exists a family of $\omega$-Hamiltonian vector fields $X:V\to V$ satisfying $dX_0=L$.
\end{corollary}

\Dem Based on the proof of Proposition \ref{relacaomatrizcampo2}, we define the symmetric matrix $B := [{\omega}]^{T}L$ and $H(x) = \frac{1}{2}\mathbf{x}^TB\mathbf{x}+F(x)$, for any function $F:V\to{\mathbb{R}}$ such that $j^2F(0)=0$. Thus, the $\omega$-Hamiltonian vector field $X$ associated with $H$ satisfies
	\[X(x) = ({[\omega]}^{-1})^{T}\nabla H(x)=({[\omega]}^{-1})^{T}B\mathbf{x}+({[\omega]}^{-1})^{T}\nabla F(x) = L\mathbf{x}+ G(x),\] where $G(x) = ({[\omega]}^{-1})^{T}\nabla F(x)$ is an $\omega$-Hamiltonian vector field such that $j^1G(0) = 0$. Hence, we have $dX_0 = L$.\cqd  

\begin{example}\label{exemplocorrecao}{\rm
		In order to illustrate the previous corollary, we consider the matrices\begin{center}
			$L=\left[
			\begin{array}{ccccccc}
				0 & 1 &  &  &  &  &  \\
				0 & 0 &  &  &  &  &  \\
				&  & 0 & \lambda_1 &  &  &  \\
				&  & -\lambda_1 & 0 &  &  &  \\
				&  &  &  & \ddots &  &  \\
				&  &  &  &  & 0 & \lambda_{n-1} \\
				&  &  &  &  & -\lambda_{n-1} & 0 \\
			\end{array}\right]$
		\end{center}and $[\omega]\in{\mathbb{M}}_{2n}({\mathbb{R}})$ as in Example \ref{exemplos Cap4}. Note that $L\in{\mathfrak{sp}}_{\omega}(n;{\mathbb{R}})$. Let $$B = [\omega]^TL = \begin{bmatrix}
			0 &  &  &  &  &  & & &   \\
			& a_1  &  &  &  &  & & &  \\
			&  &  & a_2\lambda_1 &  &  & & &  \\
			&  &  &  &  a_2\lambda_1  &  & & &  \\
			&  &  &  & & \ddots &  & & \\
			&  &  &  &  &  & & a_n\lambda_{n-1} & \\
			&  &  &  &  & & & &  a_n\lambda_{n-1} \\
		\end{bmatrix}$$		
		and define $H:{\mathbb{R}}^{2n}\to{\mathbb{R}}$ by $H(x) = \dfrac{1}{2}\mathbf{x}^TB\mathbf{x}+F(x)$, where $F(x_1, \ldots, x_{2n}) = x_1^2 \sin x_1.$ By the proof of Corollary \ref{relacaomatrizcampo2barra}, we have that 
		\begin{align*}
			X(x) & = L\mathbf{x}+([\omega]^{-1})^T\nabla F(x) \\ & = L\mathbf{x} + (0, -\dfrac{1}{a_1}(2x_1 \sin x_1 + x_1^2 \cos x_1), 0 , \ldots, 0) \\ & = (x_2, -\dfrac{1}{a_1}(2x_1 \sin x_1 + x_1^2 \cos x_1), \lambda_1 x_4, - \lambda_1 x_3, \ldots, \lambda_{n-1}x_{2n}, - \lambda_{n-1}x_{2n -1}) 
		\end{align*} is an $\omega$-Hamiltonian vector field associated with $H$, with $dX_0=L$.
}\end{example}

\quad

We prove now that to construct an $\omega$-Hamiltonian vector field $X$ satisfying the conditions of Corollary \ref{relacaomatrizcampo2barra}, it is enough to consider a skew-symmetric matrix $L\in GL(2n)$. In another direction, a natural question that arises is: given an arbitrary vector field $X: V \to V$, 
is there a symplectic structure $\omega$ such that $X$ is $\omega$-Hamiltonian? The next corollary and Proposition \ref{relacaomatrizcampo2} provide us a partial answer to this question in the linear context.

\begin{corollary}\label{matrizantissimetrica}	
	Let $L\in GL(2n)$ be a skew-symmetric matrix. Then $L\in{\mathfrak{sp}}_{\omega}(n;{\mathbb{R}})$, for the bilinear form $\omega$ such that $[\omega] = -L^{-1}$. 
\end{corollary}

\Dem Consider the symplectic vector space $(V,\omega)$ where $[\omega] = -L^{-1}$ and the symmetric matrix $B = I_{2n}$. Proposition \ref{relacaomatrizcampo1} states that there exists an $\omega$-Hamiltonian vector field $X$ such that $$dX_0 =([\omega]^{-1})^TI_{2n} =  -L^T = L.$$ Thus, $L\in{\mathfrak{sp}}_{\omega}(n;{\mathbb{R}})$, by Proposition \ref{relacaomatrizcampo3}. \cqd  

We finish this section with some results relative to polynomial Hamiltonian vector fields. We denote by $\mathcal{P}^{k}$ the vector space of all homogeneous  polynomial functions $V\to\mathbb{R}$ of degree $k$ and by $\vec{\mathcal{P}}^k$ the vector space of all homogeneous polynomial maps $V\to V$ of degree $k$. A monomial in $\vec{\mathcal{P}}^k$ is an expression of the form $x^{\alpha}\mathbf{v}_j$, where  $\mathbf{v}_j$ is the $j$-th element of a basis of $V$. If $\vec{\mathcal{P}}$ denotes the vector space of all polynomial maps $V\to V$, then $$\vec{\mathcal{P}} = \displaystyle\bigoplus_{k=0}^{\infty}\vec{\mathcal{P}}^k.$$ 

\begin{lemma}\label{dXxsimetrica} (\cite[Chapter 2, Lemma 7.5]{Chow}) For each $k\geq2$, given $f\in\vec{\mathcal{P}}^k$ such that $df_x$ is a symmetric matrix for all $x\in V$, there exists a unique function $H^{k+1}\in\mathcal{P}^{k+1}$ satisfying\[\nabla H^{k+1}(x) = f(x).\]\end{lemma}

The next theorem establishes a relation between polynomial $\omega$-Hamiltonian vector fields $X$ and matrices $dX_x$. However, it is important to note that the necessary condition for $X:V\to V$ to be an $\omega$-Hamiltonian vector field is true even if $X$ is not polynomial (see Proposition \ref{relacaomatrizcampo3}).

\begin{theorem}\label{relacaopolinomiais}Let $X\in\vec{\mathcal{P}}$ such that $X(0) = 0$. Then $X$ is an $\omega$-Hamiltonian vector field if and only if $dX_{x}\in{\mathfrak{sp}}_{\omega}(n;{\mathbb{R}})$ for all $x\in V$.
\end{theorem}

\Dem The necessary condition follows by Proposition \ref{relacaomatrizcampo3}. Conversely, suppose that $dX_{x}\in{\mathfrak{sp}}_{\omega}(n;{\mathbb{R}})$ for all $x\in V$. Since $X(0) = 0$, we write $X = X^1+X^2+\cdots+X^k$, where $X^j\in\vec{\mathcal{P}}^j$, $j=1,\ldots,k$. Thus $dX_x = dX^1_x+dX^2_x+\cdots+dX^k_x$, whence\begin{center}
		$\begin{array}{rcl}
			0 &= & (dX_{x})^T[\omega]+[\omega] dX_x\\
			&=&\left(dX_{x}^1+dX_{x}^2+\cdots+dX_{x}^k\right)^T[\omega]+[\omega] \left(dX_{x}^1+dX_{x}^2+\cdots+dX_{x}^k\right)\\
			&=&\left[(dX_{x}^1)^T[\omega]+[\omega] dX_x^1\right]+\cdots+\left[(dX_{x}^k)^T[\omega]+[\omega] dX_x^k\right].
		\end{array}$
	\end{center}Each sum in square brackets is a matrix, where each entry is zero or a homogeneous polynomial of degree $j-1$, $j=1,\ldots,k$. So 
	\begin{equation}
		\label{eq:bracket}
		(dX_{x}^j)^T[\omega]+[\omega] dX_x^j=0,
	\end{equation} for each $j=1,\ldots,k$. Thus $dX_{x}^j\in{\mathfrak{sp}}_{\omega}(n;{\mathbb{R}})$. In particular, $X^1 = dX_{0}^1\in{\mathfrak{sp}}_{\omega}(n;{\mathbb{R}})$. According to Proposition \ref{relacaomatrizcampo2}, $X^1$ is an $\omega$-Hamiltonian vector field, that is, there exists $H^2\in{\mathcal{P}}^2$ such that $$X^1(x) =X_{H^2}(x) = ([\omega]^{-1})^{T}\nabla H^{2}(x).$$ For each $j\geq2$, it follows from (\ref{eq:bracket}) that the matrix $[\omega]^TdX_x^j$ is symmetric. By Lemma \ref{dXxsimetrica} there exists $H^{j+1} \in \mathcal{P}^{j+1}$ such that $\nabla H^{j+1}(x) = [\omega]^TX ^j(x)$. If \begin{center}
		$H(x) = \displaystyle\sum_{j=1}^{k} H^{j+1}(x)$,
	\end{center} then $X(x) = ([\omega]^{-1})^{T}\nabla H(x)$ for all $x\in V$. Hence, $X$ is an $\omega$-Hamiltonian vector field.\cqd  

\begin{example}\label{naohamiltonianamatriznula}{\rm
		Consider the polynomial vector field\[X(x_1,x_2,x_3,x_4) = (x_1^2x_3+x_1x_2x_4, x_1x_2x_3+x_2^2x_4, x_1x_3^2+x_2x_3x_4, x_1x_3x_4+x_2x_4^2).\] We have that $X$ is not $\omega$-Hamiltonian, although its linearization $dX_0$ is an $\omega$-Hamiltonian matrix, for any symplectic structure $\omega$ on ${\mathbb{R}}^4$. Indeed, note that
		\begin{center}
			$dX_x = \left[\begin{array}{cccc}
				2x_1x_3+x_2x_4&x_1x_4&x_1^2&x_1x_2\\
				x_2x_3&2x_2x_4+x_1x_3&x_1x_2&x_2^2\\
				x_3^2&x_3x_4&2x_1x_3+x_2x_4&x_2x_3\\
				x_3x_4&x_4^2&x_1x_4&2x_2x_4+x_1x_3
			\end{array}\right]$.
		\end{center} Since $dX_0$ is the null matrix, then $dX_0\in{\mathfrak{sp}}_{\omega}(n;{\mathbb{R}})$ for all $\omega$. However, for $x = (1,-1,1,2)$ the matrix
		\begin{center}
			$dX_x = \left[\begin{array}{cccc}
				0&2&1&-1\\
				-1&-3&-1&1\\
				1&2&0&-1\\
				2&4&2&-3
			\end{array}\right]$
		\end{center}is not $\omega$-Hamiltonian, because its trace is non-zero. By Theorem \ref{relacaopolinomiais}, $X$ is not an $\omega$-Hamiltonian vector field for any $\omega$.
}\end{example}


\section*{Acknowledgments} 

Eralcilene Moreira Terezio was financed in part by the Coordena\c{c}\~ao de Aperfei\c{c}oamento de Pessoal de N\'ivel Superior - Brasil (CAPES) - Finance 001.


\begin{thebibliography}{99}
	\bibitem{Alomair} Alomair R., Montaldi  J. (2017). Periodic orbits in Hamiltonian systems with involutory symmetries. \emph{J. Dyn. Differ. Equ.} 29 (4), 1283-1307. 
	
	\bibitem{Terezio} Baptistelli P. H., Rodrigues Hernandes M. E., Terezio E. M. Normal forms of $\omega$-Hamiltonian vector fields with symmetries. In preparation.
	
	\bibitem{buono} Buono P.-L., Laurent-Polz F., Montaldi J. (2005). Symmetric Hamiltonian bifurcations. In: Montaldi J., Ratiu T., ed(s). \emph{Geometric Mechanics and Symmetry: The Peyresq Lectures}. London Mathematical Society Lecture Note Series. Vol. 306. Cambridge, UK: Cambridge University Press, 357-402.
	
	\bibitem{buzzi1} Buzzi C. A., Teixeira M. A. (2004). Time-reversible Hamiltonian vector fields with symplectic symmetries. \emph{J. Dyn. Differ. Equ.} 16 (2) 559-574. 
	
	\bibitem{Chow} Chow S.-N., Li C., Wang D. (1994). \emph{Normal forms and bifurcation of planar vector fields}. Cambridge, UK: Cambridge University Press.
		
	\bibitem{Hofer} Hofer H., Zehnder E. (1994). \emph{Symplectic Invariants and Hamiltonian Dynamics}. Birkh\"{a}user Advanced Texts Basler Lehrb\"{u}cher. Berlin, DE: Birkh\"{a}user, Basel.
	
	\bibitem{NJ} Jacobson N. (1985). \emph{Basic algebra I}. 2nd ed. Mineola, US: Dover Publications.
	
	\bibitem{lee} Lee J. M. (2013). \emph{Introduction to Smooth Manifolds}. Graduate Texts in Mathematics. Vol. 218, 2nd ed. New York, US: Springer-Verlag.
	
	\bibitem{Ricardo2013} Martins R. M. (2013). Formal equivalence between normal forms of reversible and Hamiltonian dynamical systems. \emph{Commun. Pure Appl. Anal.} 13(2), 703-713. 
	
	\bibitem{Ricardo2011} Martins R. M., Teixeira M. A. (2011). On the similarity of Hamiltonian and reversible vector fields in 4D. \emph{Commun. Pure Appl. Anal.} 10 (4), 1257-1266. 
	
	\bibitem{McDuff} McDuff D., Salamon D. (2017). \emph{Introduction to Symplectic Topology}. Oxford Graduate Texts in Mathematics. Vol. 27, 3rd ed. Oxford, UK: Oxford University Press.
	
	\bibitem{Montaldi2000} Montaldi J., Roberts M. (2000). Note on semisymplectic actions of Lie groups. \emph{Compt. Rendus Acad. Sci. Math.} 330 (12), 1079-1084. 
	
	\bibitem{Montaldi} Montaldi J., Roberts M., Stewart I. (1990). Existence of nonlinear normal modes of symmetric hamiltonian systems. \emph{Nonlinearity.} 3 (3), 695-730. 
	
	\bibitem{Robinson} Robinson R. C. (1971). \emph{Lectures on Hamiltonian systems}. Monografias de Matem\'atica. Vol. 7. Guanabara, BR: IMPA.
	
	\bibitem{Sev} Sevryuk M. B. (1986). \emph{Reversible Systems}. Lectures Notes in Mathematics. Vol. 1211. Berlin, DE: Springer-Verlag.
	
	\bibitem{Cannas} Silva, A. C. (2008). \emph{Lectures on symplectic geometry}. Lecture Notes in Mathematics. Vol. 1764. Berlin, DE: Springer-Verlag.
	
	\bibitem{Claudia} Wulff C., Roberts M. (2002). Hamiltonian systems near relative periodic orbits. \emph{SIAM J. Appl. Dyn. Syst.} 1 (1), 1-43. 
	
	\bibitem{Zehnder} Zehnder E. (2010). \emph{Lectures on dynamical systems: Hamiltonian vector fields and symplectic capacities}. EMS Textbooks in Mathematics. Vol 11. Freiburg, DE: European Mathematical Society.
	
\end{thebibliography}
\end{document}